\documentclass{amsart}
\usepackage{amsmath, amssymb, amscd, mathtools, graphicx, amsfonts, enumerate, mathrsfs}

\newtheorem{theorem}{Theorem}[section]
\newtheorem{proposition}[theorem]{Proposition}
\newtheorem{lemma}[theorem]{Lemma}

\theoremstyle{definition}

\newtheorem{example}[theorem]{Example}
\newtheorem{remark}[theorem]{Remark}

\numberwithin{equation}{section}

\newcommand{\N}{\mathbb{N}}                        
\newcommand{\Z}{\mathbb{Z}}                        
\newcommand{\K}{\mathbb{K}}                        
\newcommand{\R}{\mathbb{R}}                        
\newcommand{\C}{\mathbb{C}}                        

\newcommand{\OO}{\mathcal{O}}                      
\newcommand{\supp}{\mathrm{supp}}                  
\newcommand{\inform}{\mathrm{in}}                  
\newcommand{\inexp}{\mathrm{exp}}                  

\newcommand{\NN}{\mathfrak{N}}                     
\newcommand{\DD}{\mathscr{D}}                      
\newcommand{\mm}{\mathfrak{m}}                     

\begin{document}

\title{On finite determinacy of complete intersection singularities}

\author{Janusz Adamus}
\address{Department of Mathematics, The University of Western Ontario, London, Ontario, Canada N6A 5B7}
\email{jadamus@uwo.ca}
\author{Aftab Patel}
\address{Department of Mathematics, The University of Western Ontario, London, Ontario, Canada N6A 5B7}
\email{apate378@uwo.ca}
\thanks{J. Adamus's research was partially supported by the Natural Sciences and Engineering Research Council of Canada}


\begin{abstract}
We give an elementary combinatorial proof of the following fact: Every real or complex analytic complete intersection germ $X$ is equisingular -- in the sense of the Hilbert-Samuel function -- with a germ of an algebraic set defined by sufficiently long truncations of the defining equations of $X$.
\end{abstract}
\maketitle

\section{Introduction}
\label{sec:intro}

The question of finite determinacy is one of the central problems in singularity theory. When dealing with singularities of (real or complex) analytic sets or mappings, one would often like to forget the original infinite transcendental data and to work instead with its (sufficiently long) Taylor truncation. This approach is satisfactory in many circumstances. For example, the Milnor number of an isolated hypersurface singularity can be correctly calculated this way. 
In general, however, local analytic invariants of a given singularity may differ from those of its Taylor approximations of arbitrary length (see Example~\ref{ex:no-equality} below).

The present paper is concerned with complete intersection singularities. It seems not so well known that, from the algebraic point of view, complete intersection singularities are finitely determined. More precisely, as in Theorem~\ref{thm:main} below, the Hilbert-Samuel function of a (real or complex) analytic germ defined by a regular sequence $\{f_1,\dots,f_k\}$ coincides with the Hilbert-Samuel function of the germ defined by sufficiently long Taylor polynomials of the series $f_1,\dots,f_k$. In this sense, every transcendental complete intersection singularity is equisingular with a germ of an algebraic set.
This result follows from the work of Srinivas and Trivedi \cite{ST}. Here, we give an elementary alternative proof.

\subsection{Main results}

Let $\K=\R$ or $\C$. Let $x=(x_1,\dots,x_m)$ and let $\mm_x$ denote the maximal ideal in the ring of convergent power series $\K\{x\}$. For a natural number $\mu\in\N$ and a power series $f\in\K\{x\}$, the \emph{$\mu$-jet} of $f$, denoted $j^\mu f$, is the image of $f$ under the canonical epimorphism $\K\{x\}\to\K\{x\}/\mm_x^{\mu+1}$.
For an ideal $J$ in $\K\{x\}$, let
\[
H_J(\eta)=\dim_\K\K\{x\}/(J+\mm_x^{\eta+1})\,,\quad\eta\in\N
\]
denote the Hilbert-Samuel function of $\K\{x\}/J$.

\begin{theorem}
\label{thm:main}
Let $X$ be a $\K$-analytic subspace of $\K^m$, of dimension $m-k$ at $0\in X$. Suppose that the local ring $\OO_{X,0}=\K\{x\}/I$ is a complete intersection, and $\{f_1,\dots,f_k\}$ is a regular sequence in $\K\{x\}$ which generates the ideal $I$. Then, there exists $\mu_0\in\N$ such that, for every $\mu\geq\mu_0$ and for every $k$-tuple $\{g_1,\dots,g_k\}\subset\K\{x\}$ satisfying $j^\mu{g_i}=j^\mu{f_i}$, $i=1,\dots,k$, we have:
\begin{itemize}
\item[(i)] The $k$-tuple $\{g_1,\dots,g_k\}$ is a regular sequence in $\K\{x\}$
\item[(ii)] The ideal $J\coloneqq(g_1,\dots,g_k)\!\cdot\!\K\{x\}$ satisfies $H_J(\eta)=H_I(\eta)$ for all $\eta\in\N$.
\end{itemize}
\end{theorem}

Our proof of Theorem~\ref{thm:main} is elementary. Our approach is combinatorial, via the so called \emph{diagrams of initial exponents} of Hironaka (see Section~\ref{sec:diagram-ini} for details). In fact, Theorem~\ref{thm:main} is a straightforward consequence of our main result, Theorem~\ref{thm:AS-conj} below, concerning stabilization of the sequence of diagrams of initial exponents of ideals $I_\mu$ which are Taylor approximations of a given ideal $I$ in $\K\{x\}$.

More precisely, given an ideal $I$ in $\K\{x\}$, generated by some power series $f_1,\dots,f_k$, one defines $I_\mu$ to be the ideal generated by the $\mu$-jets $j^\mu{f_1},\dots,j^\mu{f_k}$. As was shown in \cite{AS3}, the diagram $\NN(I)$ of initial exponents of the ideal $I$ is then contained in the diagram $\NN(I_\mu)$ of $I_\mu$, for all $\mu$ sufficiently large. Since the $I_\mu$ are generated by polynomials, it is desirable to know if there exists $\mu$ large enough so that $\NN(I_\mu)=\NN(I)$. Theorem~\ref{thm:AS-conj} asserts that this is indeed the case when $f_1,\dots,f_k$ form a regular sequence. This gives an affirmative answer to a recent conjecture of Adamus-Seyedinejad (\cite[Conj.\,3.7]{AS3}).

\subsection{Plan of the paper}

As mentioned above, our main tool here is Hironaka's diagram of initial exponents. We recall this notion and its relevance to the Hilbert-Samuel function in Section~\ref{sec:diagram-ini}. Simply speaking, calculating the Hilbert-Samuel function of a quotient $\K\{x\}/I$ amounts to counting the points in the complement of the diagram of $I$ (cf. Remark~\ref{rem:HS-function}).

To make this work easily accessible to a wide audience, we recall in Section~\ref{sec:prelim} the basic notions from local algebra and analytic geometry used in the paper. We also show there how the problem stated in Theorem~\ref{thm:main} over a field $\K$, which is either $\R$ or $\C$, always reduces to the complex case.

Section~\ref{sec:complement} contains the key combinatorial argument of the paper, Proposition~\ref{prop:key}. Combined with Proposition~\ref{prop:HS-mult} below, it allows one to relate the multiplicity of the ring $\K\{x\}/I$ to the cardinality of the so-called generic level of the complement of the diagram of $I$.

Section~\ref{sec:approx} is concerned with approximation of diagrams. The main results, Theorems~\ref{thm:AS-conj} and~\ref{thm:main}, are proved in Section~\ref{sec:proofs}.

\section{Preliminaries}
\label{sec:prelim}

A sequence of elements $a_1,\dots,a_k$ in a ring $A$ is called a \emph{regular sequence} on $A$ if the ideal $(a_1,\dots,a_k)$ is proper, $a_1$ is a non-zerodivisor in $A$ and, for each $i=1,\dots,k-1$, the image of $a_{i+1}$ is a non-zerodivisor in $A/(a_1,\dots,a_i)$.
Recall (\cite[Ch.\,VIII, \S\,9, Cor.\,2]{ZS}) that if $A$ is a local ring and $a\in A$ is a non-zerodivisor then $\dim A/(a)=\dim A-1$ (where $\dim$ denotes the Krull dimension).

A ring $R$ is called a \emph{complete intersection} if there is a regular ring $A$ and a regular sequence $a_1,\dots,a_k$ in $A$ such that $R\cong A/(a_1,\dots,a_k)$. In particular, if $I=(a_1,\dots,a_k)$ is an ideal in a regular local ring $A$ of dimension $m$, then $A/I$ is a complete intersection if and only if its Krull dimension satisfies \,$\dim A/I=m-k$.

We have the following (see, e.g., \cite[Ch.\,VIII, \S\,8--10]{ZS} or \cite[\S\,13--14]{M}):

\begin{proposition}
\label{prop:HS-mult}
For an ideal $I$ in a Noetherian local ring $A$, the Hilbert-Samuel function $H_I(\eta)$ of $A/I$, for sufficiently large $\eta\in\N$, is a polynomial of degree $d=\dim A/I$ \,in $\eta$, whose initial coefficient is of the form $e(I)/d!$, where $e(I)\in\Z$. The integer $e(I)$ is called the \emph{multiplicity} of the ring $A/J$.
\end{proposition}

Let now $\K=\R$ or $\C$, let $x=(x_1,\dots,x_m)$, and let $\K\{x\}$ denote the ring of convergent power series in variables $x$ with coefficients in $\K$.
Let $\{f_1,\dots,f_k\}$ be a regular sequence in $\K\{x\}$, and let $I\coloneqq(f_1,\dots,f_k)\cdot\K\{x\}$.
Let $X$ be a complex analytic subspace of $\C^m$ whose local ring at $0\in\C^m$ is defined by the ideal $I\cdot\C\{x\}$; i.e., $\OO_{X,0}\cong\C\{x\}/I\!\cdot\!\C\{x\}$. By the Macaulay unmixedness theorem (see, e.g., \cite[Cor.\,18.14]{Eis}), all associated primes of $I$ in $\C\{x\}$ are of height $k$, and hence all irreducible components of the germ $X_0$ are of dimension $m-k$. In particular, the germ $X_0$ is reduced, and so $X$ can be thought of simply as a complex analytic subset of an open neighbourhood $\Omega$ of $0$ in $\C^n$, which is of pure dimension $m-k$.

Then, after a linear change of coordinates in $\C^m$, there is a fundamental system of neighbourhoods $U=U'\times U''$ of $0\in\C^m$, with $U'\subset\C^{m-k}$ and $U''\subset\C^k$, such that the restriction $\pi|_X:X\cap U\to U'$ of the projection $\pi:U'\times U''\to U'$ is a proper and surjective map, and $(\pi|_X)^{-1}(0)=\{0\}$ (see, e.g., \cite[Ch.\,III, Prop.\,4]{Nar}).
Let $p$ denote the cardinality of a generic fibre of $\pi|_X$. By \cite[Thm.\,6.5]{D}, we have
\begin{equation}
\label{eq:Draper}
p=e(I)\,,
\end{equation}
where $e(I)$ is the multiplicity of the local ring $\OO_{X,0}\cong\C\{x\}/I\!\cdot\!\C\{x\}$.

Since $X$ is pure-dimensional and the fibres of $\pi|_X$ are finite, the Remmert open mapping theorem (see, e.g., \cite[Ch.\,V, \S\,6, Thm.\,2]{Loj}) implies that $\pi|_X$ is open. It then follows from the Cohen-Macaulayness of $\OO_{X,0}$ and \cite[Prop.\,3.20]{F} that $\pi|_X$ is a flat mapping (after shrinking $U$, if needed). Finally, recall that, by \cite[Cor.\,3.13]{F}, a finite complex analytic map $\varphi:X\to Y$, with $Y$ reduced,  is flat if and only if the multiplicity map
\[
\nu_\varphi:Y\ni y\mapsto\nu_\varphi(y)=\!\!\sum_{x\in\varphi^{-1}(y)}\!\!\dim_\C\OO_{\varphi^{-1}(y),x}\in\Z
\]
is locally constant on $Y$. Since, over a generic $y\in U'$, $\nu_{\pi|_X}(y)$ is just the cardinality $p$ of the fibre $(\pi|_X)^{-1}(y)$, it follows from \eqref{eq:Draper} that
\begin{multline}
\label{eq:mult1}
e(I)=\dim_\C\OO_{(\pi|_X)^{-1}(0),0}=\dim_\C\C\{x_1,\dots,x_k\}/(I\!\cdot\!\C\{x\})(0)=\\
\dim_\K\K\{x_1,\dots,x_k\}/I(0)\,,
\end{multline}
where the evaluation is at $x_{k+1}=\dots=x_m=0$ (cf. Section~\ref{sec:diagram-ini}). The last equality in \eqref{eq:mult1} follows from the fact that $(I\!\cdot\!\C\{x\})(0)=I(0)\!\cdot\!\C\{x_1,\dots,x_k\}$.

\section{Diagram of initial exponents and Hilbert-Samuel function}
\label{sec:diagram-ini}

In this section, we recall the notion of Hironaka's diagram of initial exponents as well as his division theorem. In fact, we shall only use it here in the following simplified setting. For a detailed exposition, we refer the reader to \cite{BM}.

Let $\K=\R$ or $\C$. Let $x=(x_1,\dots,x_m)$ and let $\mm_x$ denote the maximal ideal in the ring of convergent power series $\K\{x\}$.
We will write $x^\beta$ for $x_1^{\beta_1}\dots x_m^{\beta_m}$, where $\beta=(\beta_1,\dots,\beta_m)\in\N^m$.

Let $1\leq k<m$ be an integer. We will sometimes distinguish the last $m-k$ variables and write $\tilde{x}=(x_{k+1},\dots,x_m)$, for short.
In that case, for a power series $F\in\K\{x\}=\K\{\tilde{x}\}\{x_1,\dots,x_k\}$, we define its \emph{evaluation at $0$} as
$F(0)=F(x_1,\dots,x_k,0,\dots,0)\in\K\{x_1,\dots,x_k\}$, and for an ideal $J$ in $\K\{x\}$
define $J(0)\coloneqq\{F(0): F\in J\}$, the \emph{evaluated} ideal.

We define a total ordering of $\N^m$ by lexicographic ordering of the $(m+1)$-tuples $(|\beta|,\beta_1,\dots,\beta_m)$, where $\beta=(\beta_1,\dots,\beta_m)$ and $|\beta|\coloneqq\beta_1+\dots+\beta_m$ is the \emph{length} of $\beta$.
The \emph{support} of $F=\sum_{\beta\in\N^m}f_\beta x^\beta$ is defined as $\supp(F)\coloneqq\{\beta\in\N^m:f_\beta\neq0\}$.
The \emph{initial exponent} of $F$, denoted $\inexp(F)$, is the minimum (with respect to the above total ordering) over all $\beta\in\supp(F)$.
Similarly, $\supp(F(0))=\{(\beta_1,\dots,\beta_m)\in\supp(F):\beta_{k+1}=\dots=\beta_m=0\}$ and $\inexp(F(0))=\min\{\beta\in\supp(F(0))\}$, for the evaluated
series (with respect to the total ordering induced on $\N^{k}$). Of course, $\supp(F(0))\subset\supp(F)$.


Given an ideal $J$ in $\K\{x\}$, we denote by $\NN(J)$ the \emph{diagram of initial exponents} of $J$, that is,
\[
\NN(J)=\{\inexp(F):F\in J\setminus\{0\}\}\,.
\]
Similarly, for the evaluated ideal $J(0)$, we set
\[
\NN(J(0))=\{\inexp(F(0)):F\in J, F(0)\neq0\}\,.
\]
Note that every diagram $\NN(J)\subset\N^m$ satisfies the equality $\NN(J)+\N^m=\NN(J)$. (Indeed, for $\beta\in\NN(J)$ and $\gamma\in\N^m$, one can choose $F\in J$ such that $\inexp(F)=\beta$; then $x^\gamma F\in J$ and hence $\beta+\gamma=\inexp(x^\gamma F)$ is in $\NN(J)$.)

\begin{remark}
\label{rem:vertices}
Let $\DD(m)=\{\NN\subset\N^m:\NN+\N^m=\NN\}$ be the collection of \emph{diagrams} in $\N^m$.
It is not difficult to show that, for every $\NN\in\DD(m)$, there exists a unique smallest (finite) set $V(\NN)\subset\NN$ such that $V(\NN)+\N^m=\NN$ (see, e.g., \cite[Lem.\,3.8]{BM}). The elements of $V(\NN)$ are called the \emph{vertices} of the diagram $\NN$.
\end{remark}

We now recall a combinatorial interpretation of Hironaka's division theorem:
For a proper ideal $J$ in $\K\{x\}$, set $\Delta=\N^m\setminus\NN(J)$, and define $\K\{x\}^\Delta=\{F\in\K\{x\}:\supp(F)\subset\Delta\}$. Consider the canonical projection $\K\{x\}\to\K\{x\}/J$ and its restriction to $\K\{x\}^\Delta$, called $\kappa$.

\begin{proposition}[{cf. \cite[\S\,6, Prop.\,9]{H}}]
\label{prop:Hir1}
The mapping $\kappa:\K\{x\}^\Delta\to\K\{x\}/J$ is surjective.
In other words, every power series \,$F\in\K\{x\}\!\setminus\! J$ is congruent modulo $J$ to a power series supported in $\N^m\setminus\NN(J)$.
\end{proposition}

\begin{remark}
\label{rem:HS-function}
The above proposition allows one to express the Hilbert-Samuel function of an ideal in terms of the complement of the diagram of initial exponents of that ideal:
Let $x=(x_1,\dots,x_m)$, and let $\mm_x$ be the maximal ideal in $\K\{x\}$.
For an ideal $J$ in $\K\{x\}$, let $H_J(\eta)=\dim_\K\K\{x\}/(J+\mm_x^{\eta+1})$ denote the Hilbert-Samuel function of $\K\{x\}/J$.
It follows from Proposition~\ref{prop:Hir1} that
\begin{equation}
\label{eq:HS-function}
H_J(\eta)\,=\,\#\,(\N^m\setminus\NN(J))\cap\{\beta\in\N^m:|\beta|\leq \eta\}.
\end{equation}
\end{remark}

We complete this section with the following simple but useful observation.

\begin{proposition}
\label{prop:diagram-dim}
For an ideal $J$ in $\K\{x\}$, the following conditions are equivalent:
\begin{itemize}
\item[(i)] $\dim(\K\{x\}/J)\leq\dim\K\{x\}-k$.
\item[(ii)] After a linear change of coordinates in $\K^n$, the diagram $\NN(J)$ has a vertex on each of the first $k$ coordinate axes in $\N^m$.
\end{itemize}
\end{proposition}

\begin{proof}
Condition (ii) clearly implies (i). On the other hand, (i) implies that (after a linear change of coordinates, if needed) $\K\{x\}/J$ is a finite $\K\{\tilde{x}\}$-module, where $\tilde{x}=(x_{k+1},\dots,x_m)$ (see, e.g., \cite[Ch.\,III, Prop.\,2]{Nar}). The latter is equivalent to saying that the images of $x_1,\dots,x_k$ in $\K\{x\}/J$ are integral over $\K\{\tilde{x}\}$. Therefore, by Proposition~\ref{prop:Hir1}, for every $i=1,\dots,k$, the complement of the diagram $\NN(J)$ in $\N^m$ contains at most finitely many elements on the $i$'th coordinate axis in $\N^m$. Hence (ii).
\end{proof}

\section{Counting points in the complement of a diagram}
\label{sec:complement}

Let $k$ and $m$ be positive integers, with $k<m$. For a diagram $\NN\in\DD(m)$, set $\Delta(\NN)\coloneqq\N^m\setminus\NN$. Define
\[
\DD_k(m)\coloneqq\{\NN\in\DD(m): \exists\alpha\in\Z_+\mathrm{\ s.\,t.\ }(\alpha,0,\dots,0),\,\dots,\,(0,\dots,\stackrel{(k)}{\alpha},0,\dots)\in\NN\}\,.
\]
Then, $\NN\in\DD_k(m)$ if and only if $\Delta(\NN)\subset\{0,\dots,\alpha-1\}^k\times\N^{m-k}$ for some $\alpha\in\Z_+$. Equivalently, $\NN$ has a vertex on each of the first $k$ coordinate axes in $\N^m$ (cf. Remark~\ref{rem:vertices}).
Further, let $\DD^*_k(m)$ denote the set of those $\NN\in\DD_k(m)$ that have no vertices on any other coordinate axis of $\N^m$.

For $\NN\in\DD^*_k(m)$ and $a=(a_{k+1},\dots,a_m)\in\N^{m-k}$, define
\[
L_a(\NN)\coloneqq\{(\beta_1,\dots,\beta_k)\in\N^k: (\beta_1,\dots,\beta_k,a_{k+1},\dots,a_m)\in\Delta(\NN)\}\,,
\]
and let $\delta_a(\NN)\coloneqq\,\#\,L_a(\NN)$ denote the cardinality of $L_a(\NN)$. We will call $L_a(\NN)$ the \emph{$a$-level} of $\Delta(\NN)$.

\begin{remark}
\label{rem:generic-level}
Note that, by finitness of the vertex set $V(\NN)$ (Remark~\ref{rem:vertices}), for every $\NN\in\DD^*_k(m)$ there exists $N\in\N$ such that
\[
L_a(\NN)=L_{a'}(\NN)\mathrm{\ \ for\ all\ } a,a'\in\N^{m-k}\setminus\{(\beta_{k+1},\dots,\beta_m)\;:\;\beta_i<N,\ i=k+1,\dots,m\}\,.
\]
We may thus speak of the \emph{generic level} $L_a(\NN)$ of $\Delta(\NN)$.
\end{remark}

The following result is the key technical ingredient of our arguments.

\begin{proposition}
\label{prop:key}
Let $k$ and $m$ be positive integers, with $k<m$.
Let $\NN\in\DD^*_k(m)$, and let $\delta$ be the cardinality of the generic level of $\Delta(\NN)$.
Then, for sufficiently large $\eta\in\N$, the function
\[
\Phi_\NN(\eta)\coloneqq\,\#\,\Delta(\NN)\cap\{\beta\in\N^m:|\beta|\leq\eta\}
\]
is a polynomial in $\eta$ of degree $m-k$ with initial coefficient $\dfrac{\delta}{(m-k)!}$\;.
\end{proposition}

For the proof of the proposition, we will need the following simple observation, which we prove for the sake of completeness.

\begin{lemma}
\label{lem:1}
Let $S_t^d$ denote the number of $d$-tuples $(\beta_1,\dots,\beta_d)\in\N^d$ satisfying
\[
\beta_1+\dots+\beta_d\leq t\,,\quad\mathrm{where\ }t\in\N\,.
\]
Then, $S_t^d$ is a polynomial of degree $d$ in $t$ with leading coefficient $1/d!$
\end{lemma}

\begin{proof}
We proceed by induction on $d$. For $d=1$, clearly $S_t^1=t+1$, as required. Suppose then that $d+1\geq2$ and we have
\[
S_t^d = \frac{1}{d!}t^d+a_{d-1}t^{d-1}+\dots+a_0\,.
\]
To find $S_t^{d+1}$, let $\tilde{S}_{\xi}$ be the number of solutions in 
$\mathbb{N}^d$ to
\[
\beta_1+\dots+\beta_d\leq t-\xi\,.
\]
Observe that $S_t^{d+1}=\sum_{\xi=0}^{t}\tilde{S}_{\xi}$.
By the inductive hypothesis, we have
\[
\tilde{S}_{\xi}=S_{t-\xi}^{d}=\frac{1}{d!}(t-\xi)^d+a_{d-1}(t-\xi)^{d-1}+\dots+a_0\,,
\]
hence, after rearranging the terms of $S_t^d,\dots,S_0^d$,
\begin{equation}
\label{eq:sum1}
S_t^{d+1} = \frac{1}{d!}\left(\sum_{i=0}^t i^d\right)+a_{d-1}\left(\sum_{i=0}^t i^{d-1}\right)+\dots+a_0\left(\sum_{i=0}^{t} 1\right)\ . 
\end{equation}

Now, for $n,\nu\geq1$, consider the sum $S_{\nu}(n)\coloneqq\sum_{i=0}^n i^{\nu}$. We shall show, by induction on $\nu$, that $S_\nu(n)$ is a polynomial in $n$ of degree $\nu+1$ with leading coefficient $1/(\nu+1)$. Indeed, for $\nu=1$, we have $S_1(n)=n(n+1)/2$, as required. For $\nu\geq2$, the identity 
\[
(p+1)^{\nu+1} - p^{\nu+1} = \sum_{r=0}^{\nu} \binom{\nu+1}{r} p^r
\]
summed up over $p=1,\dots,n$, yields
\[
(n + 1)^{\nu+1}-1 = \sum_{r=0}^{\nu}\binom{\nu+1}{r}\sum_{p = 1}^n p^r = \sum_{r=0}^{\nu}\binom{\nu+1}{r}S_r(n)\,,
\]
hence
\begin{equation}
\label{eq:sum2}
S_{\nu}(n) = \frac{(n+1)^{\nu+1}-1}{\nu+1} - \frac{\sum_{r=0}^{\nu-1}\binom{\nu+1}{r}S_r(n)}{\nu+1}\,.
\end{equation}
By the inductive hypothesis, the second summand on the right hand side of \eqref{eq:sum2} is a polynomial in $n$ of degree $\nu$, and hence the degree and leading coefficient of $S_\nu(n)$, as determined by the first summand, are $\nu+1$ and $1/(\nu+1)$, respectively.

Finally, applying the above to \eqref{eq:sum1} with $\nu=t$, we obtain
\[
S_t^{d+1} = \frac{1}{d!}\cdot\frac{1}{d+1}\;t^{d+1} + \phi(t)\,,
\]
where $\phi$ is a polynomial in $t$ of degree less than $d+1$.
\end{proof}
\medskip

We are now ready to prove Proposition~\ref{prop:key}.

\begin{proof}[Proof of Proposition~\ref{prop:key}]
Let $N\in\N$ be such that $L_a=L_{a'}$ for all $a,a'\in\N^{m-k}\setminus\Gamma$, where $\Gamma=\{(\beta_{k+1},\dots,\beta_m)\;:\;\beta_i<N,\ i=k+1,\dots,m\}$ (see Remark~\ref{rem:generic-level}). Pick $a\in\N^{m-k}\setminus\Gamma$.

By finiteness of $\Gamma$, there is a constant $C$ such that $\Phi_\NN(\eta)=C+\sum_{i=1}^{m-k}P_i(\eta)$, where $P_1(\eta)$ is the number of $m$-tuples $(\beta_1,\dots,\beta_m)$ in $\N^m$ satisfying
\begin{eqnarray*}
\beta_1+\dots+\beta_m &\leq& \eta\\ \beta_{k+1} &\geq& N\\ (\beta_1,\dots,\beta_k) &\in& L_a\,,
\end{eqnarray*}
and, for $i>1$, $P_i(\eta)$ is the number of $m$-tuples $(\beta_1,\dots,\beta_m)$ in $\N^m$ satisfying
\begin{eqnarray*}
\beta_1+\dots+\beta_m &\leq& \eta\\
\beta_{k+1} &<& N\\
\vdots \\
\beta_{k+i-1} &<& N\\
\beta_{k+i} &\geq& N\\
(\beta_1, \dots, \beta_k) &\in& L_a\,.
\end{eqnarray*}
It now suffices to show that, for $\eta$ sufficiently large, $P_1(\eta)$ is a polynomial in $\eta$ of degree $m-k$ with initial coefficient $\delta_a/(m-k)!$, and each $P_i(\eta)$, for $i>1$, is a polynomial in $\eta$ of degree strictly less than $m-k$.

First, let us consider $P_1(\eta)$. By applying a coordinate transformation $\beta_{k+1}\coloneqq\beta_{k+1}+N$, we see that $P_1(\eta)$ is the same as the number of $m$-tuples satisfying
\begin{eqnarray*}
\beta_1 + \dots + \beta_m &\leq& \eta - N \\
(\beta_1, \dots, \beta_k) &\in& L_a.
\end{eqnarray*}
We define $C_\nu$ to be the number of $m$-tuples $(\beta_1,\dots,\beta_m)$ in $\N^m$ satisfying
\begin{eqnarray*}
\beta_1 + \dots + \beta_m &\leq& \eta - N\\
\beta_1 + \dots + \beta_k &=& \nu\\
(\beta_1, \dots, \beta_k) &\in& L_a.
\end{eqnarray*}
Further, let $B_\nu$ be the number of $(m-k)$-tuples $(\beta_{k+1},\dots,\beta_m)$ in $\N^{m-k}$ satisfying
\[
\beta_{k+1}+\dots+\beta_m\ \leq\ \eta-\nu-N,
\]
and let $D_\nu$ be the number of $k$-tuples $(\beta_1,\dots,\beta_k)$ in $L_a$ satisfying
\[
\beta_1+\dots+\beta_k\ =\ \nu.
\]
Then, for every $\nu\in\N$, we have $C_\nu=B_\nu D_\nu$. Note also that $P_1(\eta)=\sum_{\nu=0}^{\infty} C_\nu$.

Now, since $\NN$ has a vertex on each of the first $k$ coordinate axes in $\N^m$, there exists $M>0$ such that $D_\nu=0$ for all $\nu\geq M$, and hence
\begin{equation}
\label{eq:delta}
\delta_a=\sum_{\nu=0}^M D_\nu\,,\quad\mathrm{and}\quad P_1(\eta)=\sum_{\nu=0}^M B_\nu D_\nu\,.
\end{equation} 
Note that, in terms of Lemma~\ref{lem:1}, $\displaystyle{B_\nu=S_{\eta-\nu-N}^{m-k}}$, and thus, by that lemma,
\[
P_1(\eta)\ =\ \sum_{\nu=0}^M D_\nu\!\cdot\!\left(\frac{(\eta-\nu-N)^{m-k}}{(m-k)!}\;+\;\phi(\eta-\nu-N)\right)\,,
\]
where $\phi$ is a polynomial of degree strictly less than $m-k$.
It follows that the leading coefficient of $P_1(\eta)$ is equal to $\dfrac{\sum_{\nu=0}^{M}D_{\nu}}{(m-k)!}$\,, which is $\dfrac{\delta_a}{(m-k)!}$\,, by \eqref{eq:delta}, as required.
\smallskip

Next, we show that $P_i(\eta)$ is a polynomial of degree less than $m-k$, for any $i>1$.
Indeed, for $i>1$, let $Q_i = \{(\alpha_1,\dots,\alpha_{i-1})\in\N^{i-1}: \alpha_j<N \mathrm{\ for\ } j=1,\dots,i-1\}$,
and let $F_\mu$ be the number of $(i-1)$-tuples $(\beta_{k+1}, \dots, \beta_{k+i-1})$ in $Q_i$ satisfying
\[
\beta_{k+1}+\dots+\beta_{k+i-1}=\mu\,.
\]
Let $R\in\N$ be such that $F_\mu=0$ for $\mu\geq R$, and let $D_\nu$ be defined as above.
Also, we apply the coordinate transformation $\beta_{k+i}\coloneqq\beta_{k+i}+N$, as 
before. Let $\bar{B}_{\nu, \mu}$ be the the number of $(m-k-i+1)$-tuples $(\beta_{k+i},\dots,\beta_m)$ in $\N^{m-k-i+1}$ satisfying
\[
\beta_{k+i}+\dots+\beta_m\leq\eta-\nu-\mu-N\,,
\]
and let $\bar{C}_{\mu, \nu}$ be the number of $m$-tuples $(\beta_1,\dots,\beta_m)$ in $\N^m$ satisfying
\begin{eqnarray*}
\beta_1+\dots+\beta_m &\leq& \eta-N\\
\beta_1+\dots+\beta_k &=& \nu\\
(\beta_1,\dots,\beta_k) &\in& L_a\\
\beta_{k+1}+\dots+\beta_{k+i-1} &=& \mu\\
(\beta_{k+1},\dots,\beta_{k+i-1}) &\in& Q_i\,.
\end{eqnarray*}
Then, $\bar{C}_{\mu,\nu} = D_{\nu}\bar{B}_{\mu,\nu}F_{\mu}$, for all $\nu,\mu\in\N$, and $P_i(\eta) = \sum_{\mu = 0}^{\infty} \sum_{\nu = 0}^{\infty} \bar{C}_{\mu, \nu}$. Thus, by the choice of $M$ and $R$,
\[
P_i(\eta) = \sum_{\nu=0}^M\sum_{\mu=0}^R D_{\nu}\bar{B}_{\nu,\mu}F_{\mu}\,.
\]
Note that, in terms of Lemma~\ref{lem:1}, $\bar{B}_{\nu,\mu} = S_{\eta-\nu-\mu-N}^{m-k-i+1}$, and hence, by that lemma and because $i>1$, we get $\deg{P_i}<m-k$, which completes the proof.
\end{proof}

\section{Approximation of diagrams}
\label{sec:approx}

Let $\K=\R$ or $\C$. Let $x=(x_1,\dots,x_m)$ and let $\mm_x$ denote the maximal ideal of $\K\{x\}$. Recall that, for a natural number $\mu\in\N$ and a power series $f\in\K\{x\}$, the \emph{$\mu$-jet} of $f$, denoted $j^\mu f$, is the image of $f$ under the canonical epimorphism $\K\{x\}\to\K\{x\}/\mm_x^{\mu+1}$.

In the present section we study the relations between the diagram of initial exponents of a given ideal in $\K\{x\}$ and those of its Taylor approximations. Throughout this section, we will use the following notation: Let $f_1,\dots,f_k$ be a finite collection of power series in $\K\{x\}$ and let
\[
I=(f_1,\dots,f_k)\cdot\K\{x\}\,.
\]
For a natural number $\mu$, let $I_\mu$ denote the ideal generated by the $\mu$-jets $j^\mu{f_i}$, $i=1,\dots,k$, that is,
\[
I_\mu=(j^\mu{f_1},\dots,j^\mu{f_k})\cdot\K\{x\}\,.
\]

The following simple observation will be used often in our considerations.

\begin{remark}
\label{rem:jet-exp}
Given a power series $F\in\K\{x\}$, suppose that $\mu\geq|\inexp(F)|$.
Then 
\[
\inexp(G)=\inexp(F)
\]
for every $G\in\K\{x\}$ with $j^\mu G=j^\mu F$.
\end{remark}

Let us recall now a results from \cite{AS3} describing the connection between the diagram of initial exponents of $I$ and those of its approximations $I_\mu$. We include a short proof for the reader's convenience.

\begin{lemma}[{cf.\,\cite[Lem.\,3.2]{AS3}}]
\label{lem:diagram-up-to-l}
Let $I$ and $\{I_\mu\}_{\mu\in\N}$ be as above. Let $l_0$ be the maximum of lengths of vertices of the diagram $\NN(I)$. Then:
\begin{itemize}
\item[(i)] For every $\mu\geq l_0$ and every $k$-tuple $\{g_1,\dots,g_k\}$ satisfying $j^\mu g_i=j^\mu f_i$, $i=1,\dots,k$, the ideal $J\coloneqq(g_1,\dots,g_k)\!\cdot\!\K\{x\}$ satisfies $\NN(J)\supset\NN(I)$.\\
In particular, $\NN(I_\mu)\supset\NN(I)$ for all $\mu\geq l_0$.
\item[(ii)] Given $l\geq l_0$, for all $\mu\geq l$ and every $k$-tuple $\{g_1,\dots,g_k\}$ satisfying $j^\mu g_i=j^\mu f_i$, $i=1,\dots,k$, the ideal $J\coloneqq(g_1,\dots,g_k)\!\cdot\!\K\{x\}$ satisfies
\[
\NN(J)\cap\{\beta\in\N^m:|\beta|\leq l\}=\NN(I)\cap\{\beta\in\N^m:|\beta|\leq l\}\,.
\]
\end{itemize}
\end{lemma}

\begin{proof}
Fix $\mu\geq l_0$ and let $g_1,\dots,g_k\in\K\{x\}$ be such that $j^\mu g_i=j^\mu f_i$, $i=1,\dots,k$. 
By Remark~\ref{rem:vertices}, for the proof of (i) it suffices to show that the vertices of $\NN(I)$ are contained in $\NN(J)$.
Let then $F\in I$ be a representative of a vertex of $\NN(I)$. We can write $F=\sum_{i=1}^kh_if_i$, for some $h_i\in\K\{x\}$. Then,
\[
j^\mu F=j^\mu(\sum_{i=1}^kh_if_i)=j^\mu(\sum_{i=1}^kh_i\!\cdot\!j^\mu{f_i})=j^\mu(\sum_{i=1}^kh_i\!\cdot\!j^\mu{g_i})
=j^\mu(\sum_{i=1}^kh_ig_i)\,,
\]
since the power series of a product up to order $\mu$ depends only on the power series up to order $\mu$ of its factors.
Hence, by Remark~\ref{rem:jet-exp}, we have equality of the initial exponents $\inexp(F)=\inexp(\sum_{i=1}^kh_ig_i)$. It follows that $\inexp(F)\in\NN(J)$, which proves (i).

For the proof of part (ii), fix $l\geq l_0$. Let $\mu\geq l$ and let $g_1,\dots,g_k\in\K\{x\}$ be such that $j^\mu g_i=j^\mu f_i$, $i=1,\dots,k$. By part (i), it now suffices to show that
\[
\NN(J)\cap\{\beta\in\N^m:|\beta|\leq l\}\ \subset\ \NN(I)\cap\{\beta\in\N^m:|\beta|\leq l\}\,.
\]
Pick $\beta^*\in\N^m\setminus\NN(I)$ with $|\beta^*|\leq l$. Suppose that $\beta^*\in\NN(J)$. Then, one can choose $G\in J$ with $\inexp(G)=\beta^*$. Write $G=\sum_{i=1}^kh_i\cdot g_i$ for some $h_i\in\K\{x\}$.
We have
\[
j^\mu G=j^\mu(\sum_{i=1}^kh_ig_i)=j^\mu(\sum_{i=1}^kh_i\cdot j^\mu{g_i})=j^\mu(\sum_{i=1}^kh_i\cdot j^\mu{f_i})=j^\mu(\sum_{i=1}^kh_if_i)\,,
\]
and since $\mu\geq l\geq|\inexp(G)|$, it follows that $\inexp(G)=\inexp(\sum_{i=1}^kh_if_i)$, by Remark~\ref{rem:jet-exp} again. Therefore $\beta^*\in\NN(I)$; a contradiction.
\end{proof}

\begin{lemma}
\label{lem:finite-complement}
Let $I=(f_1,\dots,f_k)\cdot\K\{x\}$ be such that the diagram $\NN(I)$ has finite complement in $\N^m$ (i.e., $\NN(I)\in\DD_m(m)$).
Then, there exists $\mu_0\in\N$ such that, for all $\mu\geq\mu_0$ and all $k$-tuples $\{g_1,\dots,g_k\}$ satisfying $j^\mu g_i=j^\mu f_i$, $i=1,\dots,k$, the ideal $J\coloneqq(g_1,\dots,g_k)\cdot\K\{x\}$ satisfies $\NN(J)=\NN(I)$.
In particular, $\NN(I_\mu)=\NN(I)$ for all $\mu\geq\mu_0$.
\end{lemma}

\begin{proof}
Let $l_0$ be the maximum of lengths of vertices of the diagram $\NN(I)$, let $l_1\coloneqq\max\{|\inexp(f_i)|:i=1,\dots,k\}$,
and let $l_2\coloneqq\max\{|\beta|:\beta\in\N^m\setminus\NN(I)\}+1$. Set $\mu_0\coloneqq\max\{l_0,l_1,l_2\}$.

Pick $\mu\geq\mu_0$ and $g_1,\dots,g_k\in\K\{x\}$, such that $j^\mu g_i=j^\mu f_i$ for $i=1,\dots,k$. Set $J\coloneqq(g_1,\dots,g_k)\cdot\K\{x\}$. Then, Remark~\ref{rem:jet-exp} and inequality $\mu\geq l_1$ imply that $\inexp(g_i)=\inexp(f_i)$ for $i=1,\dots,k$.

Let $F\in I$ be a representative of a vertex $\beta^*$ of $\NN(I)$. Then, $F=\sum_{i=1}^kh_if_i$, for some $h_1,\dots,h_k\in\K\{x\}$. Since $\mu\geq l_0$, we have $|\inexp(F)|\leq\mu$ and hence, by Remark~\ref{rem:jet-exp}, $\inexp(F)=\inexp(j^\mu(F))$. Therefore,
\begin{multline}
\notag
\inexp(F)=\inexp(j^\mu(\sum_{i=1}^kh_if_i))=\inexp(j^\mu(\sum_{i=1}^kh_ij^\mu f_i))=\\
\inexp(j^\mu(\sum_{i=1}^kh_ij^\mu g_i))=\inexp(j^\mu(\sum_{i=1}^kh_ig_i))\,,
\end{multline}
and thus $\inexp(F)=\inexp(\sum_{i=1}^kh_ig_i)$, by Remark~\ref{rem:jet-exp} again. It follows that $\beta^*\in\NN(J)$, and hence $\NN(I)\subset\NN(J)$, since $\beta^*$ was an arbitrary vertex.

Conversely, let $G\in J$ be a representative of a vertex $\tilde\beta$ of $\NN(J)$. Then, $G=\sum_{i=1}^kh_ig_i$, for some $h_1,\dots,h_k\in\K\{x\}$. The inequality $|\inexp(G)|\leq\mu_0$ now follows from the definition of $\mu_0$ and the inclusion $\N^m\setminus\NN(J)\subset\N^m\setminus\NN(I)$ proved above. One shows as above that then $\inexp(G)=\inexp(\sum_{i=1}^kh_if_i)$, by Remark~\ref{rem:jet-exp}, and hence $\tilde\beta\in\NN(I)$. Since $\tilde\beta$ was an arbitrary vertex, we get $\NN(J)\subset\NN(I)$, which completes the proof.
\end{proof}

Note that, in general, there need not be equality between the diagrams of $I$ and $I_\mu$, for $\mu$ arbitrarily large. This is shown in Example~\ref{ex:no-equality} below. In \cite{AS3}, the authors conjectured that the equality holds for large $\mu$ in case when $I$ is a complete intersection. This is indeed the case. More generally, we have the following result.

\begin{theorem}
\label{thm:AS-conj}
Suppose that $I$ is an ideal in $\K\{x\}$ generated by a regular sequence $\{f_1,\dots,f_k\}$. Then, there exists $\mu_0\in\N$ such that, for every $\mu\geq\mu_0$ and for every $k$-tuple $\{g_1,\dots,g_k\}$ in $\K\{x\}$ satisfying $j^\mu{g_i}=j^\mu{f_i}$, $i=1,\dots,k$, we have:
\begin{itemize}
\item[(i)] The $k$-tuple $\{g_1,\dots,g_k\}$ forms a regular sequence in $\K\{x\}$
\item[(ii)] After a linear change of coordinates in $\K^m$ which makes $\K\{x\}/I$ into a finite $\K\{\tilde{x}\}$-module, the ideal $J\coloneqq(g_1,\dots,g_k)\!\cdot\K\{x\}$ satisfies $\NN(J)=\NN(I)$.
\end{itemize}
\end{theorem}

We shall prove Theorem~\ref{thm:AS-conj} in Section~\ref{sec:proofs}.
\medskip

\begin{example}[{\cite[Ex.\,3.5]{AS3}}]
\label{ex:no-equality}
Let $I$ be an ideal in $\K\{x,y\}$ generated by $f_1$ and $f_2$ of the form
\begin{align*}
f_1 &= x^3y+xy^4+xy^5+xy^6+\dots,\\
f_2 &= x^2y^3+y^6+y^7+y^8+\dots.
\end{align*}
Then, for every $\mu\geq5$, we have $y^2\cdot j^\mu f_1-x\cdot j^\mu f_2=xy^{\mu+1}$, hence $(1,\mu+1)\in\NN(I_\mu)$. However, $(1,k)\notin\NN(I)$ for any $k\geq 1$.

We prove the latter by contradiction. Suppose there exists $F\in I$ with $\inexp(F)=(1,k_0)$ for some $k_0\in\N$. Choose $h_1,h_2\in\K\{x,y\}$ so that $F=h_1f_1+h_2f_2$. Let $ax^{\alpha_1}y^{\alpha_2}$ and $bx^{\beta_1}y^{\beta_2}$ be the initial terms of $h_1$ and $h_2$ respectively. Clearly, $\inform(h_1)\cdot\inform(f_1)+\inform(h_2)\cdot\inform(f_2)=0$, for otherwise the $x$-component of $\inexp(h_1f_1+h_2f_2)$ would not be $1$. Therefore, $ax^{\alpha_1+3}y^{\alpha_2+1}+bx^{\beta_1+2}y^{\beta_2+3}=0$. It follows that $\alpha_1+1=\beta_1$, $\alpha_2=\beta_2+2$, and $a+b=0$. Consequently,
\begin{equation}
\label{eq:zero}
\inform(h_1)\cdot f_1+\inform(h_2)\cdot f_2=0.
\end{equation}
Now, set $h_i^{(1)}\coloneqq h_i-\inform(h_i)$, $i=1,2$. By \eqref{eq:zero}, we get $h_1^{(1)}f_1+h_2^{(1)}f_2=F$. Hence, by repeating the above argument, $\inform(h_1^{(1)})\cdot f_1+\inform(h_2^{(1)})\cdot f_2=0$. We can thus set $h_i^{(2)}\coloneqq h_i^{(1)}-\inform(h_i^{(1)})$, $i=1,2$, and again obtain $h_1^{(2)}f_1+h_2^{(2)}f_2=F$.
By induction, if $h_i^{(j)}=h_i^{(j-1)}-\inform(h_i^{(j-1)})$, $i=1,2$, then
\begin{equation}
\label{eq:j}
h_1^{(j)}f_1+h_2^{(j)}f_2=F,\mathrm{\ for\ all\ } j.
\end{equation}
Note that, for every $j\geq1$, the initial exponent of $h_i^{(j+1)}$ is strictly greater than that of $h_i^{(j)}$, by construction.
Therefore, by the Krull Intersection Theorem, the sequences $(h_1^{(j)})_{j\geq1}$ and $(h_2^{(j)})_{j\geq1}$ converge to zero in the Krull topology of $\K\{x,y\}$. It follows from \eqref{eq:j} that $0\!\cdot\!f_1+0\!\cdot\!f_2=F$, hence $F=0$, which contradicts the choice of $F$.\qed
\end{example}

\section{Proofs of the main results}
\label{sec:proofs}

\begin{lemma}
\label{lem:I(0)}
Let $\{f_1,\dots,f_k\}$ be a regular sequence in $\K\{x\}$ and let $I=(f_1,\dots,f_k)\cdot\K\{x\}$. Then, there exists a positive integer $\mu_0$ such that, after a linear change of coordinates in $\K^m$, for every $\mu\geq\mu_0$ and a $k$-tuple $\{g_1,\dots,g_k\}$ in $\K\{x\}$ satisfying $j^\mu g_i=j^\mu f_i$, $i=1,\dots,k$, we have:
\begin{itemize}
\item[(i)] The $g_1,\dots,g_k$ form a regular sequence in $\K\{x\}$
\item[(ii)] The ideal $J\coloneqq(g_1,\dots,g_k)\!\cdot\K\{x\}$ satisfies $\NN(J(0))=\NN(I(0))$, where the evaluation is at $x_{k+1}=\dots=x_m=0$.
\end{itemize}
\end{lemma}

\begin{proof}
By assumption on $f_1,\dots,f_k$, we have $\dim\K\{x\}/I=m-k$.
Hence, by Proposition~\ref{prop:diagram-dim}, after a linear change of coordinates in $\K^m$, we may assume that the diagram $\NN(I)$ has a vertex on each of the first $k$ coordinate axes of $\N^m$. It follows that the complement $\N^k\setminus\NN(I(0))$ is a finite set. Note that the ideal $I(0)$ is generated by the $f_1(0),\dots,f_k(0)$.

Let now $\mu_0\in\N$ be the constant from Lemma~\ref{lem:finite-complement} (for the ideal $I(0)$). Pick $\mu\geq\mu_0$ and $g_1,\dots,g_k$ in $\K\{x\}$ such that $j^\mu g_i=j^\mu f_i$, $i=1,\dots,k$, and set $J\coloneqq(g_1,\dots,g_k)\!\cdot\K\{x\}$. We then have $j^\mu(g_i(0))=(j^\mu g_i)(0)=(j^\mu f_i)(0)=j^\mu(f_i(0))$, for $i=1,\dots,k$, and hence, by Lemma~\ref{lem:finite-complement}, the ideal $J(0)$ satisfies $\NN(J(0))=\NN(I(0))$. This proves (ii). 
Moreover, the last equality implies that the complement $\N^k\setminus\NN(J(0))$ is finite, and so the Krull dimension of $\K\{x_1,\dots,x_k\}/J(0)$ is zero. Hence, $\dim\K\{x\}/J=m-k$, which means that the $k$ generators $g_1,\dots,g_k$ of $J$ form a regular sequence in $\K\{x\}$.
\end{proof}

We are now ready to prove Theorem~\ref{thm:AS-conj}.

\begin{proof}[Proof of Theorem~\ref{thm:AS-conj}]
By assumption on $f_1,\dots,f_k$, we have $\dim\K\{x\}/I=m-k$.
Hence, by Proposition~\ref{prop:diagram-dim}, after a linear change of coordinates in $\K^m$, we may assume that $\NN(I)\in\DD^*_k(m)$.
By Proposition~\ref{prop:HS-mult}, Remark~\ref{rem:HS-function}, and Proposition~\ref{prop:key}, we thus have
\begin{equation}
\label{eq:mult2}
e(I)=\delta(\NN(I))\,,
\end{equation}
where $e(I)$ is the multiplicity of the ring $\K\{x\}/I$, and $\delta(\NN(I))$ is the cardinality of the generic level of $\Delta(\NN(I))$.
By Proposition~\ref{prop:Hir1}, $\dim_\K\K\{x_1,\dots,x_k\}/I(0)=\,\#\,(\N^k\setminus\NN(I(0)))$, and hence, \eqref{eq:mult1} and \eqref{eq:mult2} imply that
\begin{equation}
\label{eq:mult3}
\delta(\NN(I)) = \,\#\,(\N^k\setminus\NN(I(0)))\,.
\end{equation}

Let now $l_0$ be the maximum of lengths of vertices of $\NN(I)$, and let $\mu_0$ be the greater of $l_0$ and the $\mu_0$ from Lemma~\ref{lem:I(0)}. Pick $\mu\geq\mu_0$, and let $\{g_1,\dots,g_k\}$ be an arbitrary $k$-tuple in $\K\{x\}$ satisfying $j^\mu g_i=j^\mu f_i$, $i=1,\dots,k$.
By Lemma~\ref{lem:I(0)}, $g_1,\dots,g_k$ form a regular sequence in $\K\{x\}$, and the ideal $J\coloneqq(g_1,\dots,g_k)\!\cdot\K\{x\}$ satisfies $\NN(J(0))=\NN(I(0))$. Thus,
\begin{equation}
\label{eq:mult4}
\#\,(\N^k\setminus\NN(J(0))) =\,\#\,(\N^k\setminus\NN(I(0)))\,,
\end{equation}
and the finiteness of the above number implies that $\NN(J)\in\DD^*_k(m)$. We may thus repeat the first part of the proof for $J$ in place of $I$, and conclude that the equality \eqref{eq:mult3} holds for $J$ as well. Hence, by \eqref{eq:mult4},
\begin{equation}
\label{eq:mult5}
\delta(\NN(J)) = \delta(\NN(I))\,,
\end{equation}
where $\delta(\NN(J))$ is the cardinality of the generic level of $\Delta(\NN(J))$.
However, by Lemma~\ref{lem:diagram-up-to-l}(i), we have $\NN(J)\supset\NN(I)$, and hence the generic level of $\Delta(\NN(J))$ is a subset of the generic level of $\Delta(\NN(I))$. Therefore, by \eqref{eq:mult5}, they must be equal. It follows that $\NN(J)=\NN(I)$, by Lemma~\ref{lem:diagram-up-to-l}(ii), which completes the proof.
\end{proof}

\begin{remark}
\label{rem:any-K}
It is perhaps useful to know that, in fact, Theorem~\ref{thm:AS-conj} holds for an arbitrary field $\K$ of characteristic zero contained in $\C$. Indeed, all the components used in the above proof hold in this general setting, since this is the case for the Weierstrass Division Theorem (see, e.g., \cite{BM}) used implicitly in Proposition~\ref{prop:diagram-dim}. Also, for any $\K$ as above and any $\mm_x$-primary ideal $J$, we have equality of dimensions of vector spaces $\dim_\K\K\{x\}/J=\dim_\C\C\{x\}/J\cdot\C\{x\}$.
\end{remark}

\begin{proof}[Proof of Theorem~\ref{thm:main}]
Theorem~\ref{thm:main} follows immediately from Theorem~\ref{thm:AS-conj}, Remark~\ref{rem:HS-function}, and the fact that the Hilbert-Samuel function of $\K\{x\}/I$ is invariant under linear coordinate changes in $\K^m$.
\end{proof}

\begin{remark}
\label{rem:hypersurface}
The proof of Theorem~\ref{thm:AS-conj} implies immediately that in the case when $X$ is a hypersurface (i.e., when $I=(f_1)$ is a principal ideal in $\K\{x\}$), the Hilbert-Samuel function $H_I(\eta)$ is uniquely determined by the multiplicity $e(I)$. More precisely, for every $\mu\geq e(I)$ and every $g_1\in\K\{x\}$ satisfying $j^\mu{g_1}=j^\mu{f_1}$, the ideal $J\coloneqq(g_1)$ satisfies $H_J(\eta)=H_I(\eta)$ for all $\eta\in\N$. Indeed, for a principal $I$, $\delta(\NN(I))$ (and hence $e(I)$, by \eqref{eq:mult2}) is equal to the cardinality of the zero level $L_0(\NN(I))$ of $\Delta(\NN(I))$, since $\NN(I)$ has only one vertex (which after a linear change of coordinates in $\K^m$ may be assumed to lie on the first coordinate axis in $\N^m$). The length $l_0$ of this vertex is then equal to $|\inexp(f_1(0))|=|\inexp(f_1)|$.
\end{remark}

\bibliographystyle{amsplain}

\end{document}